# Exploring Gender Differences in Graduation Proficiency in Mathematics Education Using a Markov Chain Model: Implications for Economic Growth in Nigeria

**Clement Onwu Iji [a] and Joshua Abah Abah [a]**

**Abstract**: This study employs an ex-post facto research design to explore the fluctuations of gender difference in academic achievement among graduating students of mathematics education. Graduation statistics for a total of 1106 graduating students of mathematics education (923 males and 183 females) from a University in North Central Nigeria were used to design a discrete-time Markov chain model for the movement of the difference (*d*) in graduating proficiency from one range of values (states) to the other. Additional goodness of fit test ($\chi^2$ = 1.731, p = 0.99924) and t-test (t = 0.4055, p = 0.6852) unveiled that *d* has stayed much the same over the 12 graduation cycles used in the study, and that whatever factors determine the difference in academic achievement between male and female graduating students of mathematics education on a graduation cycle basis have remained much the same over the years. Further analysis of the model predicted the closure of the observed gender gap in the next 15 graduation cycles ($p^{15}$). The results of this study has specifically highlighted the fact that female graduates of mathematics education are as proficient as their male counterparts in driving value added services in and beyond the education sub-sector of the Nigerian economy. Based on the findings of this study, it was recommended that future work may consider an in-depth investigation of the sensitivity of parameters that may have influenced specific probabilities given in the model.

**Key-words**: Mathematics Education, Graduation Proficiency, Markov Model, Gender Difference, Economic Growth.

---

[a] University of Agriculture, Makurdi (Nigeria), Department of Science Edcation. Correspondence: Clement Onwu Iji. ijiclements07@yahoo.com

## 1. Introduction

Mathematics education as a field of study is concerned with the tools, methods and approaches that facilitate the practice of teaching and learning mathematics. Mathematics education, particularly at the higher education level, prepares students for quantitative and symbolic reasoning and advanced mathematical skills through general education, services, major and graduate programmes. Odili (2012) argued that the mathematics educator is concerned with curriculum development, instructional development and the pedagogy of mathematics. Mathematics education basically prepares students to become innovative mathematics instructors, professionally prepared to communicate mathematics to learners at all levels.

Concern for the role of mathematics education in the overall well-being of nations was at its climax at the formulation of the Education for All (EFA) initiative at the onset of the millennium. The EFA efforts considered mathematics education as a human right, a means of fostering creativity and change, that is, propelling learners into the unknown (Haddad & Draxler, 2002). In this regard, the National Universities Commission – NUC (2007) affirmed that mathematics education programme must seek:

i. The acquisition, development and inculcation of the proper value-orientation for the survival of the individual and the society.
ii. The development of the intellectual capacities of individuals to understand and appreciate their environments.
iii. The acquisition of both physical and intellectual skills which will enable individuals to develop into useful members of the community.
iv. The acquisition of an objective view of local and external environments.

On the basis of these articulations, Nigerian universities are to produce mathematics teachers with the knowledge, skills and attitudes which will enable them to contribute to the growth and development of their communities in particular and their nation in general. Mathematics education programmes are to encourage the spirit of enquiry, creativity and entrepreneurship in teachers, and to enhance the skills of teachers in contributing to the economic growth of the country (Abah, 2016). The extent of attainment of these targets is measured by the proficiency of graduates a mathematics education programme is turning out in each cycle.

Graduation proficiency refers to the quality of lifestyle, skills and content knowledge that higher education graduates exit the educational system with. A proficient graduate of mathematics education has passion for scholarship, originality and the ability to make creative decisions. The hallmarks of proficiency are intellectual curiosity, ability to work independently and in teams, and a high level of personal motivation. These qualities are manifested when graduates possess the kind of professional and lifelong learning skills that they require to be successful in their jobs (Pitan, 2016). This is to say that they are able to bring their education to bear in their daily experiences (Francis, 2015). Proficiency, therefore, translates in acceptable character and learning.

Presently, the major indicator of character and learning at graduation is the Cumulative Grade Point Average (CGPA). Although, many may argue that proficiency in mathematics education is too complex to assert a fixed system of measurement (Grigorenko *et al.*, 2009), the graduation CGPA has come to stay as the classificatory factor of how good a graduate is. Research has shown that more than half of the time (57.4 %), graduation CGPA-based class of degree is used as the minimum entry criteria for recruitment in Nigeria (Adedeji & Oyebade, 2016). Evidently, graduating class of degree has become the single cause for worry for employers, prospective graduates and other stakeholders, particularly in mathematics education.

Consequently, the new vision for mathematics education recognizes graduation proficiency as key to achieving full employment and poverty eradication, with focus on access, equity and inclusion, quality, and learning outcomes. In terms of equity and access, not only has there been reports of disparity in enrollment statistics based on gender, cases of wide gaps in academic attainment abounds. Mukoro (2014) observed that there seems not to be comparative increase in enrollment rate of females as for the males over the years. This trend was attested to by Oloyede and Lawal (2008) who added that male graduates continued to exceed those of females at both undergraduate and post-graduate levels. With respect to academic performance, Oloyede and Lawal reported lower percentages for females in the distribution of best graduating students and award winners. This observable gender gap in achievement has also been given mixed treatment by some researchers. For instance, Afuwape and Oludipe (2008) after examining gender differences in achievement among Integrated Science students (126 male and 127 female) reported that gender gap in academic attainment could be disappearing. The point of certainty in reports of studies of achievement with respect to gender, as gleaned from available literature, is the presence of fluctuations (Lindberg, Hyde, Petersen & Linn, 2010).

One of the surest ways to explain fluctuations in life situations, such as those observed in the gender gap in academic achievement, is to model the observable randomness as a Markov process. In the Markov process, named after A. A. Markov who began the study of this type of process in 1907, the outcome of a given experiment can affect the outcome of the next process, setting up a Markov chain. To specify a Markov chain, we define a set of states, $S = \{s_1, s_2, …, s_r\}$. The process starts in one of these states and move successively from one state to another. Each move is called a *step*. If the chain is currently in state $s_i$, then it moves to state $s_j$ at the next step with a probability denoted by $p_{ij}$, and this probability does not depend upon which states the chain was in before the current state. The type of Markov chain being considered is called the *discrete-time Markov chain*.

Mathematically, a discrete-time Markov chain is a family of random variables, $\{X_n\}, n \in \{0,1,2,...\}$ which satisfy the Markov property (Helbert, 2015). In other words,

$$Prob\{X_{n+1} = j | X_0 = x_0, ..., X_{n-1} = x_{n-1}, X_n = i\} = Prob\{X_{n+1} = j | X_n = i\} \quad (1)$$

The discrete-time Markov chain was chosen for this study is suited for modelling the fluctuation in academic achievement based on gender because transitioning of this occurrence happens at discrete times, namely at each graduation cycle. Helbert (2015) observed that some of the advantages of using discrete-time Markov chain model are that they are relatively easy to derive from successive data, it does not require deep insight into the reasons for the change but can give insight into the process, and the results from a Markov chain model are easily interpretable.

Based on (1), the probabilities $p_{ij}$ are called transition probabilities. The transition probabilities $\{p_{ij}\}$ form the transition probability matrix, *P*:

$$\begin{pmatrix} p_{00} & p_{01} & p_{02} & ... & ... \\ p_{10} & p_{11} & p_{12} & ... & ... \\ p_{20} & p_{21} & p_{22} & ... & ... \\ ... & ... & ... & ... & ... \\ ... & ... & ... & ... & ... \end{pmatrix}$$

The $\{p_{ij}\}$ have the properties:

$$p_{ij} \geq 0, \quad \text{all } i,j$$
$$\text{and } \sum p_{ij} = 1, \quad \text{all } i,j$$

## 2. Objetives of the Study

This study seek to explore gender differences in exiting CGPAs of graduating students of Mathematics Education using a discrete-time Markov chain model. Specifically, the study seek to:
  i. Present mean differences (*d*) in graduation proficiencies between male and female graduating students of Mathematics Education.
  ii. Define the finite states for the movement of the difference (*d*) in graduation proficiency.

iii. Generate the discrete-time Markov chain model for the movement of *d* over the graduation cycles along with the transition probabilities.
iv. Predict, based on the model, the graduation cycle in which the gap in achievement between male and female graduating students of Mathematics Education is expected to close.

### 3. Methods

This study adopts an ex-post facto research design. The ex-post facto design was considered appropriate for the study due to its scope of coverage in explaining existing relationships and developing trends. The data for this study come from the graduation statistics of exiting students of Mathematics Education programmes at a University in North Central Nigeria. The data covers a period of twelve (12) distinct graduation cycles. A graduation cycle occurs at the end of a semester in which there are students available for graduation. As such, out of the academic sessions considered in this study, four (4) contain single graduation cycles (i.e. only at the end of the second semester) while the other four (4) contain two graduation cycles (i.e. at the end of both first and second semesters). Exiting CGPAs of a total of 1106 graduating students, comprising 923 males and 183 females who successfully graduated between 2007 and 2015 from the Mathematics Education programmes offered by the university, were used in this study.

### 4. Results

The extracted data on difference in graduation proficiency between male and female graduating students of Mathematics Education is presented in Table 1. The graduation cycles covered the period between 2007/2008 academic session and 2014/2015 academic session. Sessions with two graduation cycles carry the letters "A" and "B", where "A" indicates the first semester graduation cycle and "B" indicates the second semester graduation cycle.

| Graduation Cycle | Difference in Graduation Proficiency (*d*) | Gender favoured in the Difference |
|---|---|---|
| 07/08 | 0.52 | Male |
| 08/09 | 0.06 | Male |
| 09/10 | 0.14 | Male |
| 10/11 | 0.06 | Male |
| 11/12/A | 0.08 | Female |
| 11/12/B | 0.36 | Female |
| 12/13/A | 0.08 | Male |
| 12/13/B | 0.04 | Male |
| 13/14/A | 0.29 | Female |
| 13/14/B | 0.10 | Female |
| 14/15/A | 0.35 | Male |
| 14/15/B | 0.01 | Female |

Table 1: *Difference in Graduation Proficiencies of Male and Female Graduating Students of Mathematics Education*

The model of interest in this study defines five (5) finite s*tates* for the movement of the difference in graduation proficiency. This difference, represented as *d*, is observed as it transitions between the classes or states in Table 2. The state $0.40 \leq d \leq 0.50$ does not exist based on the available data set.

| State | Range |
|---|---|
| s1 | $0.01 \leq d \leq 0.10$ |
| s2 | $0.10 \leq d \leq 0.20$ |
| s3 | $0.20 \leq d \leq 0.30$ |
| s4 | $0.30 \leq d \leq 0.40$ |
| s5 | $0.50 \leq d \leq 0.60$ |

Table 2: *Transition States and d Classes*

A chi-square goodness of fit for *d* in Table 1 yields $\chi^2$ value of 1.731 (p = 0.99924) which is greater than the critical $\chi^2$ value ($\chi^2_{0.05,11}$ = 19.675). This implies that from the set of observed differences in graduation proficiency, there is insufficient evidence to reject the null hypothesis that *d* is not any different from an equal difference in academic achievement based on gender each graduation cycle. There is a considerable chance of Type II error here, the error of failing to reject the null hypothesis when it is not correct. It is obvious that the difference in graduation CGPA is not equal for each of these 12 graduation cycles and as a result, the null hypothesis is not exactly correct. What the chi-square goodness of fit test shows is that there is *little difference* between the fluctuations in the difference in academic achievement according to gender and an assumption that there was no change in this difference (*d*) each graduation cycle. What might be concluded then is that *d* has stayed much the same over these 12 graduation cycles and that whatever factors determine the difference in academic achievement between male and female graduating students on a graduation cycle basis have remained much the same over these years.

Despite the outcome of the goodness of fit test, the model of interest is still relevant on the ground that a certain degree of random fluctuation in *d* over a period of time is bound to take place even if the basic factors determining this difference does not change. In the traditions of Moody and DuCloux (2014), this stochastic process of counting the movement of *d* results in the relative frequency of times $d$ began in a particular state and transitioned to each of the other states. Therefore, the transition probability can be defined by

$$p_{ij} = \frac{x_{ij}}{\sum_{n=1}^{k} x_{ij}} \qquad (3)$$

Literally, (3) is obtained as the number of times $d$ in state $i$ transitioned to state $j$, divided by the total number of times d was in state $i$. For simplification purpose, the total number of times $d$ was in each of the states can be observed as displayed in Table 3.

| State | Items (d) in state | Total number of times in state |
|---|---|---|
| s1 | 0.01, 0.04, 0.06, 0.06, 0.08, 0.08 | 6 |
| s2 | 0.10, 0.14 | 2 |
| s3 | 0.29 | 1 |
| s4 | 0.35, 0.36 | 2 |
| s5 | 0.52 | 1 |

Table 3: *Items in each state*

Using Tables 1 and 3 side by side, it is quite easy to see that $d$ transitioned from $s_1$ back to $s_1$ twice (i.e. 10/11 to 11/12/A and 12/12/A to 12/13/B), yielding a $p_{11}$ of $\frac{2}{6}$. Similarly, $d$ transitioned from $s_1$ to $s_3$ once (12/13/B to 13/14/A), giving rise to the transition probability $p_{13} = \frac{1}{6}$. Continuing with this pattern across the twelve (12) graduation cycles for the five (5) defined states will yield the transition matrix $P$.

$$P = \begin{pmatrix} 0.333 & 0.167 & 0.167 & 0.333 & 0 \\ 0.5 & 0 & 0 & 0.5 & 0 \\ 0 & 1 & 0 & 0 & 0 \\ 1 & 0 & 0 & 0 & 0 \\ 1 & 0 & 0 & 0 & 0 \end{pmatrix} \qquad (4)$$

From $P$, it can be observed that the probability of $d$ transitioning from the range 0.10 – 0.20 to 0.01 – 0.50 is 0.5. This implies that half of the time, a difference in CGPA between male and female graduating students of Mathematics Education in the region 0.10 – 0.20 is set to reduce to the region of 0.01 – 0.10.

Having established the transition matrix $P$ for the discrete-time Markov chain model of the gender difference in graduation proficiency among graduating students of Mathematics Education, it is possible to determine the behaviour of the chain after $n$ steps. In this type of model, the $ij$th entry of $P_{ij}^{(n)}$ of the matrix $P^n$ gives the probability that the Markov chain starting in state $s_i$, will be in state $s_j$ after $n$ steps. The implication of this property is that eventually, the successive probabilities stabilize or reach an equilibrium state and converge over time (Moody & DuCloux, 2014). Also, if a steady state can be reached for $P$ after a particular power of $P$, it can be interpreted that the existing difference ($d$) will close after the particular number of steps. In this model under consideration, a power of $P$ resulting in equilibrium will indicate the graduation cycle in which the gap in achievement between male and female graduating students of Mathematics Education is expected to close.

An application software, Microsoft Mathematics (© 2010 Microsoft Corporation, Version: 4.0.1108.0000), was used to generate successive powers of $P$ in quest for equilibrium. Around the 10$^{th}$ power, the transition probabilities tend to converge, with little variations. However, at $P^{15}$, equilibrium was attained for $P$. $P^{15}$ correct to four decimal places is given in (5). This steady-state matrix has all its rows the same.

$$P^{15} = \begin{pmatrix} 0.4997 & 0.1669 & 0.0835 & 0.2499 & 0 \\ 0.4997 & 0.1669 & 0.0835 & 0.2499 & 0 \\ 0.4997 & 0.1669 & 0.0835 & 0.2499 & 0 \\ 0.4997 & 0.1669 & 0.0835 & 0.2499 & 0 \\ 0.4997 & 0.1669 & 0.0835 & 0.2499 & 0 \end{pmatrix} \quad (5)$$

From (5), it is clear that after fifteen (15) graduation cycles, the predictions about $d$ are independent of the value of $d$ at the close of the 2014/2015 academic session, pointing to a close in the gap in academic achievement between male and female graduating students of Mathematics Education. Consequently, this information translates to a closure in gender difference in achievement within the next eight (8) years at the rate of two graduation cycles per year.

To verify this result, an independent samples t-test (equal variance assumed) was carried out on the raw data (graduation CGPA) of the 923 male and 183 female graduating students of Mathematics Education. The outcome of the t-test was not significant at 0.05 level of significance, as shown in Table 4.

| Gender | n | Df | t-stat | t-crit | p - value |
|---|---|---|---|---|---|
| Male | 923 | 1104 | 0.4055 | 1.9621 | 0.6852 |
| Female | 183 | | | | |

$\alpha = 0.05$

Table 4: *T-test of Exiting CGPA of Male and Female Graduating Students of Mathematics Education*

The results in Table 4 indicate that generally, the observed difference in academic performance between male and female graduating students of Mathematics Education is statistically not significant (p = 0.6852), affirming the possibility of near or total disappearance of such difference in the future.

## 5. Discussion and Implications for Economic Growth in Nigeria

Gender is a complex dynamic force that affects every social interaction, including interactions in educational and job settings. Gender inequality in learning mathematics and science has continued to be a topical issue of global concern, with no clear cut answer as to academic disparities between male and female students (Halpern *et al.*, 2007). This disparity translates into glaring difference in consideration within the labour market. Gender division in labour refers to the organization of labour on the assumption that men perform specific roles such as those of providers and breadwinners in the productive or wage labour sector outside the home, and that women provide domestic labour as housewives within the home (Association of African Universities, 2006).

Women made up a little over half the world's population but their contribution to measured economic activity, growth and well-being is far below its potential, with serious macroeconomic consequences. Female labour force participation has remained lower than male participation, women account for most unpaid work, and when women are employed in paid work, they are overrepresented in the informal sector and among the poor (Elborgh-Woytek *et al.*, 2013). This reality is more glaring in the field of Mathematics Education in which existing cultures around the world and particularly in Nigeria tend to create a situation of repeated priming of mathematics as negatively stereotyped on female students. This biased orientation, coupled with perceived difference in academic ability has dampened the chances of women in being employed in the field of practice and even outside the field.

However, Spelke (2005) held that in-depth studies yield little support to some stereotypical assertions, emphasizing that highly talented male and female students show equal abilities to learn mathematics. In the same vein, Iji, Abah and Anyor (2017) towed the line of Lindberg *et al.* (2010) in establishing that due to cultural shifts initiated by increasing levels of technology penetration in recent years, the gender gap is closing. This line of discovery is strongly supported by the results of this present study.

Current realities thrown up by the Nigerian economy has brought to the fore the call for sober reflection on a wide range of issues, including gender. To many, the dynamics playing out in the economy are the end products of unchecked profligacy and prolonged downgrading of cultural orientation (Abah, 2016). The prejudiced perspectives of employers of labour with respect to gender must change if the country must tap into the potentials of graduates of mathematics education. Apart from the skills of running a school plant, the training of the Mathematics Education graduate imbibes relevant qualities such as punctuality, honesty, hard work, smartness and innovation.

The strength of the graduate mathematics educationist lies in vital areas of communication, information technology, critical thinking, leadership, team working, problem solving and entrepreneurship. These attributes that are well sought after in government parastatals, private companies and firms across the country will be greatly missed if gender biases are permitted to pervade hiring decisions in the disguise of academic proficiency. Young female graduates of mathematics education should be considered as work-ready individuals who are prepared to take decisions, act according to instruction, find opportunities, take initiative and produce results in their work places (Adedeji & Oyebade, 2016). This study has buttress the fact that both male and female graduates of Mathematics Education must be treated fairly in the way they are adjudged competent for entrusted positions and contribution to economic development in Nigeria.

The outcome of this study underscores the possibility of closure of gender gaps in achievement among graduating students of Mathematics Education. A direct implication of this result is the realization of the provision of full human capacity for growth in all sectors of the economy. This study has specifically highlighted the existing equity in academic proficiency of graduates of Mathematics Education, underscoring the reality that female graduates are as proficient as their male counterparts in driving value-added services in and beyond the education sub-sector.

## 6. Conclusion

This study has attempted to explore the dynamism of difference in academic achievement of male and female graduating students of Mathematics Education. This was achieved via a discrete-time Markov chain model which attained stability after 15 graduation cycles, predicting a closure of the gender gap in the nearest future. The implications of this reality were discussed in line with the prospects for economic growth in Nigeria.

The data used for this study covers only twelve (12) graduation cycles and as such may be an inherent limitation of this work. Also, the study vaguely highlighted the presence of some constant factors moderating the difference in academic achievement of male and female graduating students of mathematics education without providing any insight on the extent of their influence.

Future considerations of this magnitude may seek graduation statistics across a wider span and location. Also, there is the need to investigate the sensitivity of parameters that may have influence specific transition probabilities given in the model.